\documentclass{amsart}
\usepackage{amssymb}

\newcommand{\IN}{\mathbb N}
\newcommand{\IR}{\mathbb R}
\newcommand{\IQ}{\mathbb Q}

\newcommand{\wdi}{\mathrm{wd}}

\newcommand{\sci}{\mathrm{sc}}
\newcommand{\Int}{\mathrm{Int}}
\newcommand{\w}{\omega}
\newcommand{\C}{\mathcal C}
\newcommand{\hl}{\mathrm{hl}}
\newcommand{\parn}{\mathrm{par}}
\newcommand{\U}{\mathcal{U}}
\newcommand{\cov}{\mathrm{cov}}
\newcommand{\add}{\mathrm{add}}
\newcommand{\cof}{\mathrm{cof}}
\newcommand{\dec}{\mathrm{dec}}
\newcommand{\clov}{\mathrm{dec}_c}
\newcommand{\F}{\mathcal{F}}

\newcommand{\cl}{\mathrm{cl}}

\newcommand{\M}{\mathcal M}
\newcommand{\A}{\mathcal A}
\newcommand{\sAE}{\sigma \mathrm{AE}}
\newcommand{\concat}{\hat{\phantom{v}}}
\newcommand{\bb}{\mathfrak b}
\newcommand{\cc}{\mathfrak c}

\newtheorem{thm}{Theorem}[section]
\newtheorem{lem}[thm]{Lemma}
\newtheorem{prop}[thm]{Proposition}
\newtheorem{cor}[thm]{Corollary}
\theoremstyle{definition}
\newtheorem{dfn}[thm]{Definition}
\newtheorem{ex}[thm]{Example}
\newtheorem{question}[thm]{Question}
\newtheorem{rem}[thm]{Remark}

\begin{document}

\title{On scatteredly continuous maps between topological spaces}
\thanks{Dedicated to Tsugunori Nogura on the occasion of his 60th birthday}
\author{Taras Banakh}
\author{Bogdan Bokalo}
\address{Instytut Matematyki, Akademia \'Swi\c etokrzyska, Kielce (Poland) and\newline
Department of Mathematics, Ivan Franko National University of Lviv (Ukraine)}
\email{tbanakh@yahoo.com, Tery.87@mail.ru}

\subjclass{54C08}
\keywords{Scatteredly continuous map, weakly discontinuous map, piecewise continuous map, $G_\delta$-measurable map, Preiss-Simon space}
\begin{abstract} A map $f:X\to Y$ between topological
spaces is defined to be {\em scatteredly continuous} if for each
subspace $A\subset X$ the restriction $f|A$ has a point of
continuity. We show that for a function $f:X\to Y$ from a perfectly paracompact hereditarily Baire Preiss-Simon space $X$ into a regular space $Y$ the scattered continuity of $f$ is equivalent to (i) the weak discontinuity (for each subset $A\subset X$ the set $D(f|A)$ of discontinuity points of $f|A$ is nowhere dense in $A$), (ii) the piecewise continuity ($X$ can be written as a countable union of closed subsets on which $f$ is continuous), (iii) the $G_\delta$-measurability (the preimage of each open set is of type $G_\delta$). Also under Martin Axiom, we construct a $G_\delta$-measurable map $f:X\to Y$ between metrizable separable spaces, which is not piecewise continuous. This answers an old question of V.Vinokurov.
\end{abstract}

\maketitle

\section{Introduction}

In this paper we study scatteredly continuous maps of topological spaces and show that under some restrictions on the domain and range the scattered continuity is equivalent to many other function properties considered in the literature (such as: piecewise continuity, $G_\delta$-measurability, weak discontinuity, etc.).

By definition, a map $f:X\to Y$ between topological spaces is {\em scatteredly continuous} if for each non-empty subspace $A\subset X$ the restriction $f|A$ has a point of continuity. Such maps were introduced in \cite{AB} and appear naturally in analysis, see \cite{BKMM}.

By its spirit the definition of a scatteredly continuous map resembles the classical definition of a pointwise discontinuous map due to R.Baire \cite{Baire}. We recall that a map $f:X\to Y$ is called {\em pointwise discontinuous} if for each non-empty {\em closed} subspace $A\subset X$ the restriction $f|A$ has a continuity point. Pointwise discontinuous maps can be characterized in many different ways. In particular, the following classical theorem of R.Baire is well-known:

 \begin{thm}[Baire] For a real-valued function $f:X\to \IR$ from a complete metric space $X$ the following conditions are equivalent:
\begin{enumerate}
\item $f$ is pointwise discontinuous;
\item $f$ is $F_\sigma$-measurable {\rm (in the sense that for any open set $U\subset Y$ the preimage $f^{-1}(U)$ is of type $F_\sigma$ in $X$);}
\item $f$ is of the first Baire class {\rm (i.e., $f$ is a pointwise limit of a sequence of continuous functions).}
\end{enumerate}
\end{thm}

A similar characterization holds also for scatteredly continuous maps:

\begin{thm}\label{int2} For a real-valued function $f:X\to \IR$ from a complete metric space $X$ the following conditions are equivalent:
\begin{enumerate}
\item $f$ is scatteredly continuous;
\item $f$ is weakly discontinuous {\em (which means that for every non-empty set $A\subset X$ the set $D(f|A)$ of discontinuity points of $f|A$ is nowhere dense in $A$);}
\item $f$ is piecewise continuous {\em (which means that $X$ has a countable closed cover $\C$ such that $f|C$ is continuous for every $C\in\C$);}
\item $f$ is $G_\delta$-measurable {\em (which means that for every open set $U\subset \IR$ the preimage $f^{-1}(U)$ is of type $G_\delta$ in $X$);}
\item $f$ is of the stable first Baire class {\em (which means that there is a sequence \mbox{$(f_n:X\to\IR)_{n\in\w}$} of continuous functions that stably converges to $f$ in the sense that for every $x\in X$ there is $n\in\w$ such that $f_m(x)=f(x)$ for all $m\ge n$).}
\end{enumerate}
\end{thm}

This characterization theorem is a combined result of investigations of many authors, in particular \cite{BM}, \cite{Kir}, \cite{CL2}, \cite{Sol}, and holds in a bit more general setting. Unfortunately all known generalizations of Theorem~\ref{int2} concern maps with metrizable range or domain. In this paper we treat the general non-metrizable case. Our main result is Theorem~\ref{8.1} which unifies the results of sections~\ref{s3}--\ref{s7} and says that the conditions (1)--(4) of Theorem~\ref{int2} are equivalent for any map $f:X\to Y$ from a perfectly paracompact hereditarily Baire Preiss-Simon space $X$ to a regular space $Y$. If, in addition, $Y\in\sAE(X)$, then all the five conditions of Theorem~\ref{int2} are equivalent. In Section~\ref{s9} we present some examples showing the necessity of the assumptions in Theorem~\ref{8.1}. In particular, we answer an old question of V.Vinokurov \cite{Vino} constructing a $G_\delta$-measurable map $f:X\to Y$ between separable metrizable spaces that fails to be piecewise continuous. 

\subsection{Terminology and notations}
 Our terminology and notation are
standard and follow the monographs \cite{Ar}, \cite{En}. A ``space'' always
means a ``topological space''. Maps between topological spaces can be discontinuous.

For a subset $A$ of a topological space $X$ by $\cl_X(A)$ or $\overline{A}$ we denote the closure of $A$ in $X$ while $\Int(A)$ stands for the interior of $A$ in $X$.
For a function $f:X\to Y$ between topological spaces by $C(f)$ and $D(f)=X\setminus C(f)$ we denote the sets of continuity and discontinuity points of $f$, respectively.

By $\IR$ and $\IQ$ we denote
the spaces of real and  rational numbers, respectively; $\w$ stands for the space of finite ordinals (= non-negative integers) endowed with the discrete topology. We shall identify cardinals with the smallest ordinals of the given size.

By $\w^{<\w}$ we denote the set of all finite sequences $s=(n_1,\dots,n_k)$ of non-negative integer numbers. For such a sequence $s=(n_1,\dots,n_k)$ and a number $n\in\w$ by 
$s\concat n=(n_1,\dots,n_k,n)$ we denote the concatenation of $s$ and $n$.

\section {Some elementary properties of scatteredly continuous maps}

We start recalling the definition of a scatteredly continuous map, the principal concept of this paper.
Then we shall prove some elementary properties of such maps and present some counterexamples.

\begin{dfn}\label{2.1} A map $f:X\to Y$ between
topological spaces is called {\em scatteredly continuous} if for
any non-empty subspace $A\subset X$ the restriction $f|A:A\to Y$ has a
point of continuity.
\end{dfn}

Scatteredly continuous maps can be considered as mapping counterparts of scattered topological spaces. Let us recall that a topological space $X$ is {\em scattered} if each non-empty subspace $A$ of $X$ has an isolated point. It is clear that a topological space $X$ is scattered if and only if any bijective map $f:X\to Y$ to a discrete space $Y$ is scatteredly continuous.
In fact, a bit more general result is true.

\begin{prop}\label{2.2} A space $X$ is scattered if and only if any bijective map $f:X\rightarrow Y$ to a scattered space $Y$ is scatteredly continuous.
\end{prop}

\begin{proof} The ``if'' part is trivial. To prove the ``only if'' part, assume that $f:X\to Y$ is a scatteredly continuous map to a scattered space $Y$. We need to prove that each non-empty subspace $A\subset X$ has an isolated point. Since $f$ is scatteredly continuous, the set $C(f|A)$ of continuity points of the restriction $f|A$ is dense in $A$. The image $f(C(f|A))\subset Y$, being scattered, contains an isolated point
$y_0$. By the continuity of $f|C(f|A)$, the point $x_0=(f|A)^{-1}(y_0)$ is
isolated in $C(f|A)$ and by the density of $C(f|A)$ in $A$, $x_0$ is
isolated in $A$.
\end{proof}

The scattered continuity of maps $f:X\to Y$ defined on spaces of countable tightness can be detected on countable subspaces of $X$.
We recall that a topological space $X$ has {\em countable tightness} if for every subset $A\subset X$ and any point $a\in\overline{A}$ there is a countable subset $B\subset A$ with $a\in\overline{B}$.

\begin{prop}\label{2.4} A map $f:X\to Y$ defined of a space $X$ of countable tightness is scatteredly continuous if and only if for every countable subspace $Q\subset X$ the restriction $f|Q$ has a continuity point.
\end{prop}

\begin{proof} The ``only if'' part of this proposition is trivial. To prove  the ``if'' part, assume that $f:X\to Y$ is not scatteredly continuous and $X$ has countable tightness. Since $f$ is not scatteredly continuous, there is a subset $D\subset X$ such that $f|D$ has no continuity point.

By induction on the tree $\w^{<\w}$ we can construct a sequence $(x_s)_{s\in\w^{<\w}}$ of points of $D$  such that for every sequence $s\in\w^{<\w}$ the following conditions hold:
\begin{itemize}
\item $x_s$ is a cluster point of the set $\{x_{s\concat n}:n\in\w\}$;
\item $f(x_s)$ is not a cluster point of $\{f(x_{s\concat n}):n\in\w\}$.
\end{itemize}

We start the inductive construction taking any point $x_\emptyset\in D$. Assuming that for some finite sequence $s\in\w^{<\w}$ the point $x_s\in D$ has been chosen, use the discontinuity of $f|D$ at $x_s$ to find a neighborhood $U$ of $f(x_s)$ such that $x_s$ is a cluster point of $D\setminus f^{-1}(U)$. Since $X$ has countable tightness there is a countable set $\{x_{s\concat n}:n\in\w\}\subset D\setminus f^{-1}(U)$ whose closure contains the point $x_s$. This completes the inductive construction.

Then the set $Q=\{x_s:s\in\w^{<\w}\}$ is countable and $f|Q$ has no continuity point.
\end{proof}

We finish this section with an example of two scatteredly continuous maps whose composition is not scatteredly continuous.

\begin{ex}\label{2.5} Let $f:\IR\to\IR_\IQ$ be the identity map from the real line equipped with the standard topology $\tau$ to the real line endowed with the topology generated by the subbase $\tau\cup\{\IQ\}$. Also let $\chi_\IQ:\IR_\IQ\to \{0,1\}$ be the characteristic function of the
set $\IQ$. It is easy to show that the maps $f:\IR\to\IR_\IQ$ and $\chi_\IQ:\IR_\IQ\to\{0,1\}$ are scatteredly continuous while their composition $\chi_\IQ\circ f:\IR\to\{0,1\}$ is everywhere discontinuous (and hence fails to be scatteredly continuous).
\end{ex}

\section{The index of scattered continuity}

The scattered continuity of a map $f:X\to Y$ can be measured by an ordinal index $\sci(f)$ defined as follows. Consider the decreasing transfinite sequence $(D_\alpha(f))_\alpha$ of subspaces of $X$ defined by recursion: $D_0(f)=X$ and $D_\alpha(f)=\bigcap_{\beta<\alpha}D(f|D_\beta(f))$ for an ordinal $\alpha>0$.
The transfinite sequence $(D_\alpha(f))_\alpha$ will be called the {\em discontinuity series} of $f$. The scattered continuity of $f$ implies that $D_\alpha(f)=\emptyset$ for some ordinal $\alpha$.
The smallest ordinal $\alpha$ with $D_\alpha(f)=\emptyset$ is called the {\em index of scattered continuity} of $f$ and is denoted by $\sci(f)$.

A decreasing transfinite sequence $(X_\alpha)_{\alpha<\beta}$ of subsets of a set $X$ will be called a {\em vanishing series} if $X_0=X$, $X_\beta=\emptyset$, $X_{\alpha+1}\ne X_\alpha$ for all $\alpha<\beta$, and $X_\alpha=\bigcap_{\gamma<\alpha}X_\gamma$ for each limit ordinal $\alpha\le\beta$. A typical example of a vanishing series is the  discontinuity series $(D_\alpha(f))_{\alpha\le\sci(f)}$ of a scatteredly continuous map $f$.

\begin{prop}\label{3.1} A map $f:X\to Y$ is scatteredly continuous if and only if there is a vanishing series $(X_\alpha)_{\alpha\le\beta}$ of subsets of $X$ such that $ X_{\alpha+1}\supset D(f|X_\alpha)$ for every ordinal $\alpha<\beta$. The smallest length $\beta$ of such a series is equal to $\sci(f)$, the index of scattered continuity of $f$.
\end{prop}

\begin{proof} The ``only if''  part follows from the definition of the discontinuity series $(D_\alpha(f))$ of $f$. To prove the ``if'' part, assume that $(X_\alpha)_{\alpha\le\beta}$ is a vanishing sequence of subsets of $X$ such that $X_{\alpha+1}\supset D(f|X_\alpha)$ for every ordinal $\alpha<\beta$. To show that $f$ is scatteredly continuous, take any non-empty subset $A\subset X$ and find the smallest ordinal $\alpha\le \beta$ such that $A\not\subset X_\alpha$. We claim that $\alpha$ is not a limit ordinal. Otherwise $A\not\subset X_\alpha=\bigcap_{\gamma<\alpha}X_\gamma$ would imply that  $A\not\subset X_\gamma$ for some $\gamma<\alpha$, which would contradict the definition of $\alpha$.

So $\alpha=\gamma+1$ for some ordinal $\gamma$. Take any point $a\in A\setminus X_{\gamma+1}$. It follows from the definition of $\alpha$ that $A\subset X_\gamma$. Since $f|X_\gamma$ is continuous at each point of the set $X_\gamma\setminus X_{\gamma+1}$, the restriction $f|A$ is continuous at $a$. This proves the scattered continuity of $f$.

Next, we show that $\sci(f)\le \beta$. This will follow as soon as we prove that
 $D_\alpha(f)\subset X_\alpha$ for all ordinals $\alpha\le\beta$.
This can be done by induction. For $\alpha=0$ we get $D_0(f)=X=X_0$.

Assume that for some ordinal $\alpha\le\beta$ we have proved that $D_\gamma(f)\subset X_\gamma$ for all $\gamma<\alpha$. We should prove that $D_\alpha(f)\subset X_\alpha$. If $\alpha$ is a limit ordinal, then $$D_\alpha(f)=\bigcap_{\gamma<\alpha}D(f|D_\gamma(f))=\bigcap_{\gamma<\alpha}D_\gamma(f)\subset\bigcap_{\gamma<\alpha}X_\gamma=X_\alpha$$
by the inductive hypothesis. If $\alpha=\gamma+1$ is a successor ordinal, then the inclusion $X_{\gamma+1}\supset D(f|X_\gamma)$ implies that $f|X_\gamma$ is continuous at points of $X_\gamma\setminus X_{\gamma+1}$. By the inductive hypothesis, $D_\gamma(f)\subset X_\gamma$ and hence $f|D_\gamma(f)$ is continuous at points of $D_\gamma(f)\setminus X_{\gamma+1}$. Consequently, $$D_\alpha(f)=D(f|D_{\gamma}(f))\subset X_{\gamma+1}=X_\alpha.$$
\end{proof}

The discontinuity series will help us to characterize scatteredly continuous maps via  well-orders.

\begin{prop}\label{3.2} A map $f:X\to Y$ between topological spaces is scatteredly continuous if and only if there exists a well-order $\prec$ on the set
  $X$ such that for any non-empty subset
  $A\subset X$ the restriction $f|_{A}:A\rightarrow Y$ is
  continuous at the point $x=\min(A,\prec)$.
\end{prop}

\begin{proof} The ``if'' part is trivial. To prove the ``only if'' part, take any scatteredly continuous map  $f:X\to Y$ and consider the discontinuity series $(D_\alpha(f))_{\alpha<\sci(f)}$ of $f$.
On each set $D_\alpha(f)\setminus D_{\alpha+1}(f)$ fix a well-order $\prec_\alpha$. The well-orders $\prec_\alpha$, $\alpha<\sci(f)$, compose a well-order $\prec$ on $X$ defined by $x\prec y$ if either $x,y\in D_\alpha(f)\setminus D_{\alpha+1}(f)$ for some $\alpha$ and $x\prec_\alpha y$ or else $x\in D_\alpha(f)$ and $y\notin D_\alpha(f)$ for some $\alpha$.

Given any non-empty subset $A\subset X$ we should check that $f|A$ is continuous at the point $a=\min(A,\prec)$. Find an ordinal $\alpha$ such that $a\in D_\alpha(f)\setminus D_{\alpha+1}(f)$ and observe that $A\subset D_\alpha(f)$. Since $D_{\alpha+1}(f)=D(f|D_\alpha(f))$, $a$ is a continuity point of $f|D_\alpha(f)$ and hence of $f|A$.
\end{proof}

\section{Scatteredly continuous and weakly discontinuous maps}\label{s3}

The pathology described in Example~\ref{2.5} cannot happen in the realm of regular spaces: in this section we shall show that the composition of scatteredly continuous maps between regular spaces is scatteredly continuous. For this we first consider compositions of scatteredly continuous and weakly discontinuous maps.

Following \cite{Vino}, we define a map $f:X\to Y$ to be {\em weakly discontinuous} if for every subspace $Z\subset X$ the set $D(f|Z)$ is nowhere dense in $Z$.\footnote{Weakly discontinuous functions are refered to as Baire$^*$ 1 functions in \cite{O'M} and Baire-one-star in \cite{Kir}.} 

\begin{prop}\label{myp1} The composition of two weakly discontinuous maps is weakly discontinuous.
\end{prop}

\begin{proof} Let $f:X\to Y$, $g:Y\to Z$ be two weakly discontinuous maps. To show that $g\circ f$ is weakly discontinuous, it suffices, given a non-empty subspace $A\subset X$ to find a non-empty open subset $U\subset A$ such that $g\circ f|U$ is continuous. The weak discontinuity of $f$ yields a non-empty open set $V\subset A$ such that $f|V$ is continuous. The weak discontinuity of $g$ yields a non-empty open subset $W\subset f(V)$ such that $g|W$ is continuous. By the continuity of $f|V$, the preimage $U=(f|V)^{-1}(W)$ is open in $V$ and hence in $A$. Finally, the continuity of the functions $f|U$ and $g|f(U)$ imply the continuity of $g\circ f|U$.
\end{proof}

\begin{prop}\label{myp2} The composition $g\circ f:X\to Z$ of a weakly discontinuous map $f:X\to Y$ and a scatteredly continuous map $g:Y\to Z$ is scatteredly continuous.
\end{prop}

\begin{proof} Given a non-empty subspace $A\subset X$ we should find a continuity point of $g\circ f|A$. The weak discontinuity of $f$ implies the existence of a non-empty open set $V\subset A$ such that $f|V$ is continuous. The scattered continuity of $g$ implies that the existence of a continuity point $y_0\in f(V)$ of the restriction $g|f(V)$.
Then any point $x_0\in (f|V)^{-1}(y_0)$ is a continuity point of $g\circ f|A$.
\end{proof}

\begin{rem}\label{4.3} In light of Proposition~\ref{myp2} it is interesting to note that the composition $g\circ f:X\to Z$ of a scatteredly continuous map $f:X\to Y$ and a weakly discontinuous map $g:Y\to Z$ need not be scatteredly continuous. A suitable pair of maps $f,g$ is given in Example~\ref{2.5} (it is easy to check that the characteristic map $\chi_\IQ:\IR_{\IQ}\to\{0,1\}$ is weakly discontinuous).
\end{rem}

In some other terms the following important result on interplay between scatteredly continuous and weakly discontinuous maps was established in \cite{AB} and \cite{BM}. We present a proof here for convenience of the reader.

\begin{thm}\label{4.4} A map $f:X\to Y$ from a topological space $X$ into a (regular) topological space $Y$ is scatteredly continuous if (and only if) $f$ is weakly discontinuous. \end{thm}

\begin{proof} The ``if'' part is trivial. To prove the ``only if'' part, assume that $f:X\to Y$ is a scatteredly continuous map into a regular space $Y$. Since the scattered continuity is preserved by taking restrictions, it suffices to check that the set $D(f)$ of discontinuity points of $f$ is nowhere dense. This is equivalent to saying that the interior of the set  $C(f)$ of continuity points of $f$ meets each non-empty open set $U\subset X$. Assuming the converse we would get $U\subset \overline{D(f)}$.
Fix a continuity point $x_0\in D(f)\cap U$ of the restriction $f|D(f)\cap U$ (it exists because of the scattered continuity of $f$). Since $x_0$ fails
to be a continuity point of $f$, there is
a neighborhood $Of(x_0)\subset Y$ of the point $f(x_0)$ such
that $f(Ox_0)\not\subset Of(x_0)$ for any neighborhood $Ox_0$ of
the point $x_0$ in \(X\). By the regularity of $Y$, find a
neighborhood $O^*f(x_0)$ of $f(x_0)$ such that
$\overline{O^*f(x_0)}\subset Of(x_0)$. By the continuity of the restriction $f|D(f)\cap U$ at $x_0$ there is a neighborhood
$O{x_0}\subset U$ of $x_0$ in $X$ such that $f(O{x_0}\bigcap
D(f))\in O^*f(x_0)$. It follows from the choice of the
neighborhoods $Of(x_0)$ and $O^*f(x_0)$ that
$f(O{x_0})\not\subset\overline{O^*f(x_0)}$ and hence there is a point $x_1\in Ox_0$ with $f(x_1)\notin O^*f(x_0)$. It follows from $f(D(f)\cap Ox_0)\subset O^*f(x_0)$ that $x_1\in C(f)$. Consequently, there is a neighborhood $Ox_1\subset Ox_0$ such that $f(Ox_1)\subset Y\setminus \overline{O^*f(x_0)}$. Since $\overline{D(f)}\supset U\supset Ox_1$, there is a point $x_2\in D(f)\cap Ox_1$. For this point we get $f(x_2)\subset f(Ox_0\cap D(f))\cap f(Ox_1)\subset O^*f(x_0)\cap (Y\setminus \overline{O^*f(x_0)})=\emptyset$, which is a contradiction.
\end{proof}

The following example shows that the regularity of $Y$ is
essential in the previous theorem.

\begin{ex}\label{4.5} The identity map $f:\IR\to\IR_\IQ$ from Example~\ref{2.5} shows that a scatteredly continuous map to a non-regular space need not be weakly discontinuous.
\end{ex}

Unifying Proposition~\ref{myp2} with Theorem~\ref{4.4} we get a promised

\begin{prop}\label{4.6} If $f:X\to Y$ and $g:Y\to Z$ are scatteredly continuous maps and the space $Y$ is regular, then the composition $g\circ f:X\to Z$ is scatteredly continuous.
\end{prop}

Like the scattered continuity, the weak discontinuity also can be measured by an ordinal index. Namely, given a weakly discontinuous map $f:X\to Y$ consider the decreasing transfinite sequence  $(\tilde D_\alpha(f))_\alpha$ of closed subsets of $X$ defined recursively: $\tilde D_0(f)=X$ and $\tilde D_\alpha(f)=\bigcap_{\beta<\alpha}\cl_X\big(D(f|\tilde D_\beta(f))\big)$. The transfinite series $(\tilde D_\alpha(f))_{\alpha}$ will be called the {\em weak discontinuity series} of $f$. The weak discontinuity of $f$ implies that $\tilde D_\alpha(f)=\emptyset$ for some ordinal $\alpha$. The smallest such an ordinal $\alpha$ is denoted by $\wdi(f)$ and is called the {\em index of weak discontinuity} of $f$.
The following characterization of the weak discontinuity can be proved by analogy with Proposition~\ref{3.1}.

\begin{prop}\label{4.7} A map $f:X\to Y$ is weakly discontinuous if and only if there is a vanishing series $(X_\alpha)_{\alpha\le\beta}$ of closed subsets of $X$ such that $ X_{\alpha+1}\supset D(f|X_\alpha)$ for every ordinal $\alpha<\beta$. The smallest length $\beta$ of such a series is equal to $\wdi(f)$, the index of weak discontinuity of $f$.
\end{prop}

 The index of weak discontinuity is related to the index of scattered continuity $\sci(f)$ and to the hereditary Lindel\"of number $\hl(X)$ of $X$ as follows.
We recall that $\hl(X)=\sup\{l(Y):Y\subset X\}$ where $l(Y)$, the {\em Lindel\"of number} of a space $Y$, is the smallest cardinal $\kappa$ such that each open cover of $Y$ has a subcover of size $\le \kappa$.

\begin{prop}\label{4.8} If $f:X\to Y$ is a weakly discontinuous map, then $\sci(f)\le\wdi(f)<\hl(X)^+$.
\end{prop}

\begin{proof} The inequality $\wdi(f)<\hl(X)^+$ follows from that fact that any strictly decreasing transfinite sequence of closed subsets of a topological space $X$ has length $<hl(X)^+$. The inequality $\sci(f)\le\wdi(f)$ follows from the inclusions $D_\alpha(f)\subset \tilde D_\alpha(f)$ that can be proved by induction on $\alpha$. Alternatively it can be derived from Proposition~\ref{3.1}.
\end{proof}

\begin{question}\label{4.9} Is $\sci(f)=\wdi(f)$ for any scatteredly continuous map $f:[0,1]\to\IR$?
\end{question}

\begin{question}\label{4.9a} Is $\sci(f)<\hl(X)^+$ for every scatteredly continuous map $f:X\to Y$?
\end{question}

Propositions~\ref{myp1} and \ref{myp2} can be quantified with help of ordinal indices $\sci(\cdot)$ and $\wdi(\cdot)$ as follows.
First we need to recall the definition of ordinal multiplication, which can be introduced by recursion: $\alpha\cdot 0=0$ and $\alpha\cdot(\beta+1)=\alpha\cdot\beta+\alpha$.

\begin{prop}\label{4.10} If $f:X\to Y$ is a weakly discontinuous map and $g:Y\to Z$ is a scatteredly continuous map, then
$\sci(g\circ f)\le\sci(g)\cdot\wdi(f)$.
\end{prop}

\begin{proof} Observe that for every ordinal $\alpha<\wdi(f)$ the function $f$ is continuous at the set $\tilde C_\alpha(f)=\tilde D_\alpha(f)\setminus \tilde D_{\alpha+1}(f)$, which is open and dense in $\tilde D_\alpha(f)$.

Given ordinals $\alpha<\wdi(f)$ and $\beta<\sci(g)$ consider the sets
$$\begin{aligned}
D_{\alpha,\beta}=&\tilde C_\alpha(f)\cap f^{-1}(D_\beta(g))\mbox{ and}\\
C_{\alpha,\beta}=&\tilde C_\alpha(f)\cap f^{-1}\big(C(g|D_\beta(g))\big)=D_{\alpha,\beta}\setminus D_{\alpha,\beta+1}
\end{aligned}
$$
and note that the restriction $g\circ f|D_{\alpha,\beta}$ is continuous at each point of $C_{\alpha,\beta}$.

On the product $\wdi(f)\times\sci(g)=\{(\alpha,\beta):\alpha<\wdi(f),\;\beta<\sci(g)\}$ consider the lexicographic ordering: $(\alpha,\beta)\le(\alpha',\beta')$ iff either $\alpha<\alpha'$ or $\alpha=\alpha'$ and $\beta\le \beta'$. Endowed with this ordering the product $\wdi(f)\times\sci(g)$ is order isomorphic to the ordinal product
$\sci(g)\cdot \wdi(f)$.  It is clear that the successor of a pair  $(\alpha,\beta)$ in $(\wdi(f)\times\sci(g),\le)$ is the pair $(\alpha,\beta+1)$.

For ordinals $\alpha<\wdi(f)$ and $\beta<\sci(g)$ let
$$X_{(\alpha,\beta)}=\cup\{C_{\alpha',\beta'}:(\alpha,\beta)\le(\alpha',\beta')\in\wdi(f)\times\sci(g)\}.$$
Thus we obtain a vanishing series $(X_{(\alpha,\beta)})_{(\alpha,\beta)}$ of subsets of $X$. Since the lexicographic ordering of $\wdi(f)\times \sci(g)$ has order type of the ordinal product $\sci(g)\cdot\wdi(f)$, the inequality $\sci(g\circ f)\le\sci(g)\cdot \wdi(f)$ will follow from Proposition~\ref{3.1} as soon as we prove that for each $(\alpha,\beta)\in \wdi(f)\times\sci(g)$ the restriction $g\circ f|X_{(\alpha,\beta)}$ is continuous at each point $x$ of $X_{(\alpha,\beta)}\setminus X_{(\alpha,\beta+1)}$.

Since  $x\in X_{(\alpha,\beta)}\setminus X_{(\alpha,\beta+1)}=C_{\alpha,\beta}\subset X\setminus \tilde D_{\alpha+1}(f)$, the continuity of $g\circ f|X_{(\alpha,\beta)}$ at $x$ will follow as soon as we prove the continuity of $g\circ f$ restricted to the set $Z_{\alpha,\beta}=X_{(\alpha,\beta)}\setminus \tilde D_{\alpha+1}(f)$ at $x$. For this observe that
$Z_{\alpha,\beta}=\bigcup_{\beta'\ge \beta}C_{\alpha,\beta'}=\tilde C_\alpha(f)\cap f^{-1}(D_\beta(g))$ and $f(x)\in f(C_{\alpha,\beta})\subset C(g|D_\beta(g))$. Now the continuity of $f$ on the set $Z_{\alpha,\beta}$ and the inclusion $f(Z_{\alpha,\beta})\subset D_\beta(g)$ imply that $g\circ f|Z_{\alpha,\beta}$ is continuous as $x$.
\end{proof}

By analogy one can prove

\begin{prop}\label{4.11} If $f:X\to Y$ and $g:Y\to Z$ are two weakly discontinuous maps then their composition $g\circ f$ is weakly discontinuous and
$\wdi(g\circ f)\le\wdi(g)\cdot\wdi(f)$.
\end{prop}

\begin{cor} If $f:X\to Y$ is scatteredly continuous (and weakly discontinuous) map and $Z\subset X$, then $\sci(f|Z)\le \sci(f)$ (and $\wdi(f|Z)\le\wdi(f)$).
\end{cor}

\section{Scatteredly continuous and piecewise continuous maps}

In \cite{Vino} V.A.Vinokurov  proved that a map $f:X\to Y$ defined on a (complete) metrizable space $X$ is weakly discontinuous (if and) only if $f$ is piecewise continuous. We recall that a map $f:X\to Y$ between topological spaces is {\em piecewise continuous} if $X$ has a countable closed cover $\C$ such that for every $C\in\C$ the restriction $f|C$  is continuous. In \cite{BKMM} the Vinokurov's results was generalized to maps defined on perfectly paracompact (hereditarily Baire) spaces:

\begin{thm}\label{5.1} A map $f:X\to Y$ from a perfectly paracompact (and hereditarily Baire) space $X$ to a topological space $Y$ is weakly discontinuous (if and) only if $f$ is piecewise continuous.
\end{thm}

We recall that a topological space $X$ is called
\begin{itemize}
\item {\em hereditarily Baire} if each closed subspace $F\subset X$ is Baire (in the sense that the intersection of a sequence of open dense subsets of $F$ is dense in $F$);
\item {\em perfectly paracompact} if $X$ is paracompact and each closed subset of $X$ is of type $G_\delta$ in $X$.
\end{itemize}
The class of perfectly paracompact spaces is quite large: besides metrizable spaces this class includes all hereditarily  Lindel\"of spaces and all stratifiable spaces, see \cite{Borges}.

In this section we generalize Theorem~\ref{5.1} and prove its quantitative version.

First note that the piecewise continuity of a map $f:X\to Y$ can be expressed in terms of the {\em closed decomposition number} $\clov(f)$ equal to the smallest size $|\C|$ of a closed cover $\C$ of $X$ such that $f|C$ is continuous for each $C\in\C$, see \cite{Sol}. Note that a map $f$ is continuous iff $\clov(f)=1$ and piecewise continuous iff $\clov(f)\le\aleph_0$.

The following easy proposition (whose proof is left to the reader) shows that the continuity covering number behaves nicely with respect to compositions.

\begin{prop} For any maps $f:X\to Y$ and $g:Y\to Z$ we get $\clov(g\circ f)\le \max\{\clov(g),\clov(f)\}$.
\end{prop}

It turns out that the closed decomposition number $\clov(f)$ of a weakly discontinuous map $f:X\to Y$ can be estimated from above by the index of weak discontinuity $\wdi(f)$ of $f$ and the large pseudocharacter $\Psi(X)$ of $X$. By definition, for a topological $T_1$-space $X$ the {\em large pseudocharacter}  $\Psi(X)$   is
 equal to the smallest cardinal $\kappa$ such that each closed subset $F\subset X$ can be written as the intersection $\cap\U$ of a family $\U$ consisting $\le\kappa$ open subsets of $X$. It is easy to see that $\Psi(X)\le\hl(X)$ for every regular space $X$.

\begin{prop}\label{5.2} If $f:X\to Y$ is a weakly discontinuous map from a $T_1$-space $X$, then $\clov(f)\le \max\{|\wdi(f)|,\Psi(X)\}$.
\end{prop}

\begin{proof} Given a weakly discontinuous map $f:X\to Y$ consider the weak discontinuity series $(\tilde D_\alpha(f))_{\alpha<\wdi(f)}$ of $f$. Let $\kappa=\max\{|\wdi(f)|,\Psi(X)\}$. It follows from $\Psi(X)\le\kappa$ that for every $\alpha<\wdi(f)$ the set $\tilde D_\alpha(f)\setminus \tilde D_{\alpha+1}(f)$ can be written as the union $\cup\C_\alpha$ of a family of closed subsets of $X$ with $|\C_\alpha|\le\kappa$. The continuity of $f$ on $\tilde D_\alpha(f)\setminus \tilde D_{\alpha+1}(f)$ implies the continuity of $f$ on each set $C\in\C_\alpha$. Then $\C=\bigcup_{\alpha<\wdi(f)}\C_\alpha$ is a closed cover of $X$ such that $|C|\le\kappa$ and for every $C\in\C$ the restriction $f|C$ is continuous. This proves the inequality $\clov(f)\le\kappa=\max\{|\wdi(f)|,\Psi(X)\}$.
\end{proof}

In fact, the index of weak discontinuity of $f$ in Proposition~\ref{5.2} can be replaced by the paracompactness number $\parn(X)$ of $X$.

We recall that a topological space $X$ is {\em paracompact} if each open cover of $X$ can be refined to a locally finite open cover. According to \cite[5.1.11]{En}, this is equivalent to saying that each open cover of $X$ can be refined to a locally finite closed cover. In light of this characterization it is natural to introduce a cardinal invariant measuring the paracompactness degree of a topological space $X$.

By the {\em paracompactness number} $\parn(X)$ of a topological space $X$ we understand the smallest cardinal $\kappa$ such that each open cover of $X$ can be refined by an closed cover $\F$ of $X$ that can be written as the union $\F=\bigcup_{\alpha<\kappa}\F_\alpha$ of $\kappa$ many locally finite families  $\F_\alpha$ of closed subsets of $X$. Hence a topological space $X$ is paracompact if and only if $\parn(X)\le 1$. It is easy to see that $\parn(X)\le l(X)$ for every regular space $X$.

\begin{thm}\label{5.3} If $f:X\to Y$ is a weakly discontinuous map, then
$\clov(f)\le\max\{\parn(X),\Psi(X)\}$.
\end{thm}

\begin{proof} Given a weakly discontinuous map $f:X\to Y$, consider the weak discontinuity series $(\tilde D_\alpha(f))_{\alpha<\wdi(f)}$ of $f$.

The proof of the inequality $\clov(f)\le\max\{\parn(X),\Psi(X)\}$  will be done by induction on $\wdi(f)$. If $\wdi(f)\le1$, then $f$ is continuous and hence $\clov(f)=1\le\max\{\parn(X),\Psi(X)\}$. Assume that for some ordinal $\alpha$ our theorem is proved for weakly discontinuous maps with index of weak discontinuity $<\alpha$. Assume that $f$ has $\wdi(f)=\alpha$ and let $\lambda=\max\{\parn(X),\Psi(X)\}$.

Consider the closed set $C_0=\bigcap_{\beta<\alpha}\tilde D_\beta(f)=\bigcap_{\beta+1<\alpha}\cl_X\big(D(f|\tilde D_\beta(f))\big)$ and note that $C_0=\tilde D_\alpha(f)$ if $\alpha$ is a limit ordinal and $C_0=\tilde D_\gamma(f)$ if $\alpha=\gamma+1$ is a successor ordinal. In any case the restriction $f|C_0$ is continuous.

Using the definition of the cardinal $\Psi(X)\le\lambda$, find a family $\mathcal W=\{W_k:k<\lambda\}$ of open sets of $X$ with  $C_0=\bigcap_{k<\lambda}W_k$.

Take any ordinal $k<\lambda$ and observe that $\parn(X\setminus W_k)\le\parn(X)\le\lambda$. Then the open cover $\{X\setminus\tilde D_\beta(f):\beta<\alpha\}$ of $X\setminus W_k$ admits a closed  refinement $\F_{k}=\bigcup_{i<\kappa}F_{k,i}$, where each family $\F_{k,i}$, $i<\lambda$, is locally finite in $X\setminus W_k$.

For every set $F\in\F_{k,i}$ find an ordinal $\beta<\alpha$ with $F\subset X\setminus \tilde D_\beta(f)$. By induction on $\gamma$ it can be shown that  $\tilde D_\gamma(f|F)\subset F\cap\tilde D_\gamma(f)$ for all ordinals $\gamma$. In particular, $\tilde D_\beta(f|F)=\emptyset$, which means that $\wdi(f|F)\le\beta<\alpha$.
By the induction hypothesis, $\clov(f|F)\le\max\{\parn(F),\Psi(F)\}\le\max\{\parn(X),\Psi(X)\}=\lambda$. So, we can find a closed cover $\{C_{k,i}^{F,j}:j<\lambda\}$ of $F$ such that $f|C_{k,i}^{F,j}$ is continuous for all $\alpha<\kappa$. It follows from the local finity of $\F_{k,i}$ that the family $\C_{k,i}^j=\{C_{k,i}^{F,j}:F\in\F_{k,i}\}$ is locally finite in $X$ and hence the union $C_{k,i}^j=\cup\C_{k,i}^j$ is a closed subset of $X$. Moreover, $f|C_{k,i}^j$ is continuous. Then $\C=\{C_0\}\cup\bigcup_{k,i,j<\lambda}\C_{k,i}^j$ is a closed cover of $X$ of size $\le\lambda$ witnessing that $\clov(f)\le\lambda$.
\end{proof}

Since $\max\{\parn(X),\Psi(X)\}\le\hl(X)$ for every regular space $X$, Theorem~\ref{5.3} implies

\begin{cor}\label{5.4} If $f:X\to Y$ is a weakly discontinuous map defined on a regular space $X$, then $\clov(f)\le\hl(X)$.
\end{cor}

\section{Maps of stable first Baire class}

Let us say that a sequence of maps $(f_n:X\to Y)_{n\in\w}$ {\em stably converges} to a map $f:X\to Y$ if for every $x\in X$ there is $n\in\w$ such that $f_m(x)=f(x)$ for all $m\ge n$.

A map $f:X\to Y$ between topological spaces is defined to be of {\em the stable first Baire class} if there is a sequence $(f_n)_{n\in\w}$ of continuous functions from $X$ to $Y$ that stably converges to $f$.

Maps of the stable first Baire class have been studied in \cite{CL1}, \cite{CL2}, 
\cite{BKMM}.

\begin{prop}\label{6.0} If $Y$ is a Hausdorff space, then each map $f:X\to Y$ of the stable first Baire class is piecewise continuous.
\end{prop}

\begin{proof}  Choose a sequence $(f_n:X\to Y)_{n\in\w}$ of continuous functions that stably converges to $f$. For every $n\in\w$ consider the closed subset
$$X_n=\{x\in X:\forall m\ge n\; f_m(x)=f_n(x)\}$$and note that $f|X_n\equiv f_n|X_n$ is continuous. The stable convergence of $(f_n)$ to $f$ implies that $X=\bigcup_{n\in\w}X_n$ which means that $f$ is piecewise continuous.
\end{proof}

The converse statement to Proposition~\ref{6.0} is more subtle and holds under some assumptions on the spaces $X$ and $Y$. For example, according to \cite{BKMM} a piecewise continuous map $f:X\to Y$ is of the stable first Baire class if $X$ is normal and $Y$ is real line. The main result of this section asserts that the same is true for any path-connected space $Y\in\sAE(X)$.

To define the class $\sAE(X)$ we first recall the notion of an
$AE(X)$-space. Namely, we say that a space $Y$ is an absolute
extensor for a space $X$ and denote this by $Y\in AE(X)$ if each
continuous map $f:A\to Y$ defined on a closed subspace $A\subset
X$ has a continuous extension $\tilde f:X\to Y$.

The classical Urysohn Lemma says that $\IR\in AE(X)$ for every normal space $X$. The Dugundji's Theorem \cite{Dug} says that each convex set $Y$ in a locally convex space is an absolute extensor for any metrizable space $X$ (more generally, for any stratifiable space $X$, see \cite{Borges}).

Generalizing the notion of an $AE(X)$-space to pairs $(Y,B)$ of spaces $B\subset Y$ we shall write $(Y,B)\in AE(X)$ if each continuous map $f:A\to B$ defined on a closed subspace $A\subset X$ has a continuous extension $\tilde f: X\to Y$. Hence $Y\in AE(X)$ if and only if $(Y,Y)\in AE(X)$.

We define a space $Y$ to belong to the class $\sAE(X)$ if $Y$ has a countable cover $\{Y_n:n\in\w\}$ by
closed $G_\delta$-sets such that $(Y,Y_n)\in AE(X)$ for all $n\in\w$. In particular, a space $Y$ is a $\sAE(X)$ if it admits a countable cover $\{Y_n:n\in\w\}$ by closed $G_\delta$-subspaces $Y_n\in AE(X)$.

Observe that the boundary $Y=\partial I^2$ of the square $I^2=[0,1]^2$ is not an absolute extensor for $X=I^2$. On the other hand, for each proper subcontinuum $B\subset\partial I^2$ the pair $(\partial I^2,B)\in AE(I^2)$. Since $\partial I^2$ can be covered by two proper subcontinua, we get $\partial I^2\in\sAE(I^2)$.

\begin{ex} The Sierpi\'nski carpet $Y$:

\centerline{\input serpinski.tex}

\noindent is an example of a Peano continuum that fails to be in $\sAE(I^2)$. The reason is that for each countable closed cover $\C$ of $Y$ the Baire Theorem yields a set $C\in\C$ with non-empty interior in $Y$. Then $C$ contains a loop that is not contractible in $Y$, which means that $(Y,C)\notin AE(I^2)$.
\end{ex}

\begin{thm}\label{6.1} Each piecewise continuous map $f:X\to Y$ from a normal space $X$ to a path-connected space $Y\in \sAE(X)$ is of the stable first Baire class.
\end{thm}

\begin{proof} Find a countable cover $\{Y_n:n\in\w\}$ of $Y$ by non-empty closed $G_\delta$-subsets such that $(Y,Y_n)\in AE(X)$ for all $n\in\w$. In every set $Y_n$ choose a point $y_n\in Y_n$. Since $Y$ is path-connected, for every $n\in\w$ there is a continuous map $\gamma_n:[n,n+1]\to Y$ such that $\gamma(n)=y_n$ and $\gamma_n(n+1)=y_{n+1}$.
The maps $\gamma_n$ compose a single continuous map $\gamma:[0,+\infty)\to Y$ such that $\gamma|[n,n+1]=\gamma_n$ for all $n\in\w$. Let $\IR_+=[0,+\infty)$ and denote by  $Gr(\gamma)=\{(\gamma(t),t)\in Y\times\IR_+:t\in\IR_+\}$ the graph of $\gamma$ in $Y\times \IR_+$.

Finally consider the subspace $Z=Y\times\w\cup Gr(\gamma)$ of the product $Y\times\IR_+$ and let $p:Z\to Y$ be the restriction of the projection $Y\times\IR_+\to Y$.
Define a map $s:Y\to Y\times\w\subset Z$ letting $s(y)=(y,n)$ if $y\in Y_n\setminus\bigcup_{i<n}Y_i$. It is clear that $p\circ s$ is the identity map of $Y$. Since each $Y_n$ is a closed $G_\delta$-set in $Y$, the map $s:Y\to Z$ is piecewise continuous.

Now we are able to prove that each piecewise continuous map $f:X\to Y$ is of the stable first Baire class. Observe that the composition $g=s\circ f:X\to Z$ is piecewise continuous, being the composition of two piecewise continuous maps. Hence we can find a countable closed cover $\C$ of $X$ such that map $g|C$ is continuous for all $C\in\C$. Since $g(X)\subset Y\times\w$, we may additionally assume that for every $C\in \C$ there is a number $k(C)\in\w$ such that $g(C)\subset Y\times\{k(C)\}$.

Let $\C=\{C_n:n\in\w\}$ be an enumeration of $\C$. For every $n\in\w$ let $X_n=\bigcup_{i\le n}C_i$. It is clear that the restriction $g|X_n:X_n\to Y\times\w$ is continuous. We claim that $g|X_n$ has a continuous extension $\tilde g_n:X\to Z$.

Let $K=\{k(C_i):i\le n\}$ and for every $k\in K$ let $F_k=(g|X_n)^{-1}(Y\times\{k\})$. Observe that $g(F_k)\subset Y_k\times\{k\}$ because $s(Y)\cap Y\times\{k\}\subset Y_k\times\{k\}$. It follows that the closed sets $F_k$, $k\in K$, are pairwise disjoint and $X_n=\bigcup_{k\in K}F_k$. The normality of $X$ allows us to find open neighborhoods $U_k$ of the closed sets $F_k$, $k\in K$, that have pairwise disjoint closures in $X$.

Since $X$ is normal, there is a continuous map $h:X\to [0,+\infty)$  such that $h(\overline U_k)=\{k\}$ for all $k\in K$.

For every $k\in K$ use the inclusions $g(F_k)\subset Y_k\times\{k\}$ and  $(Y,Y_k)\in AE(X)$ to find a continuous map $g_k:X\to Y\times\{k\}$ such that $g_k|F_k=g|F_k$ and $g_k(X\setminus U_k)=\{(y_k,k)\}$.

It is easy to see that the map $\tilde g_n:X\to Z$ defined by $$
\tilde g(x)=\begin{cases} g_k(x)& \mbox{if $x\in U_k$ for some $k\in K$,}\\
\big(\gamma(h(x)),h(x)\big)&\mbox{otherwise,}
\end{cases}
$$
is continuous and so is the map $\tilde f_n=p\circ \tilde g_n:X\to Y$. Since $\tilde f_n|X_n=f|X_n$ and $(X_k)$ is an increasing sequence of subsets with $X=\bigcup_{k\in\w}X_k$, the sequence $(\tilde f_n)_{n\in\w}$ stably converges to $f$ witnessing that $f$ is of the stable first Baire class.
\end{proof}

The path-connectedness is essential in Theorem~\ref{7.1}.

\begin{ex} Consider the characteristic map $f:[0,1]\to\{0,1\}$ of the half-interval $[0,\frac12)$. It is piecewise continuous but is not of the stable first Baire class (because continuous maps from $[0,1]$ to $\{0,1\}$ all are constant). Yet, $\{0,1\}\in\sAE(X)$ for any space $X$.
\end{ex}

The condition $Y\in\sAE(X)$ also is essential in Theorem~\ref{7.1}.

\begin{ex} Let $Y$ be the Sierpi\'nski carpet
in the square $X=I^2$, $y_0\in Y$ be any point and $f:X\to Y$ be the piecewise continuous function defined by $f(x)=x$ if $x\in Y$ and $f(x)=y_0$ if $x\notin Y$. We claim that this function fails to be of the stable first Baire class. Assuming the converse, find a sequence $(f_n)_{n\in\w}$ of continuous functions that stably converges to $f$. For every $n\in\w$ consider the closed set $F_k=\{x\in X:\forall m\ge n\; f_m(x)=f_n(x)\}$ and note that $\bigcup_{k\in\w}F_k=X$. The Baire Theorem implies that for some $k\in\w$ the set $F_k\cap Y$ has non-empty interior in $Y$. Then there is a map $g:\partial I^2\to F_k\cap Y$ that cannot be extended to a continuous map $\bar g:I^2\to Y$. On the other hand, the map $g$ has a continuous extension $\tilde g:I^2\to X$. Then $f_n\circ \tilde g:I^2\to Y$ is a continuous extension of $g=f\circ g=f_n\circ \tilde g|\partial I^2$ which contradicts the choice of the map $g$.
\end{ex}

\section{Scatteredly continuous and $G_\delta$-measurable maps}\label{s7}

In this section we reveal the relation of scatteredly continuous maps to $G_\delta$-measurable maps.

We define a map $f:X\to Y$ to be {\em $G_\delta$-measurable} if for every open set $U\subset Y$ the preimage $f^{-1}(U)$ is of type $G_\delta$ in $X$ (this is equivalent to saying that for every closed subset $F\subset Y$ the preimage $f^{-1}(F)$ is of type $F_\sigma$ in $X$). 

The following proposition is immediate.

\begin{prop}\label{7.1} A piecewise continuous map $f:X\to Y$ is $G_\delta$-measurable. 
\end{prop}

According to \cite{Vino} or \cite{Kir},  each $G_\delta$-measurable map $f:X\to Y$ from a complete metric space $X$ to a regular space $Y$ is weakly discontinuous.
In Theorem~\ref{7.3} below we shall show that this is still true for maps from hereditarily Baire spaces $X$ with the Preiss-Simon property.

We define a topological space $X$ to be {\em Preiss-Simon} at a point $x\in X$ if for any subset $A\subset X$ with $x\in \overline{A}$ there is a sequence $(U_n)_{n\in\w}$ of non-empty open subsets of $A$ that converges to $x$ in the sense that each neighborhood of $x$ contains all but finitely many sets $U_n$. By $PS(X)$ we denote the set of points $x\in X$ at which $X$ is Preiss-Simon. A topological space $X$ is called a {\em Preiss-Simon} space if $PS(X)=X$ (that is $X$ is Preiss-Simon at each point $x\in X$).

It is clear that each first countable space is Preiss-Simon and each Preiss-Simon space is Fr\'echet-Urysohn. A less trivial fact due to D.Preiss and P.Simon \cite{PS} asserts that each Eberlein compact space is Preiss-Simon.

The proof of Theorem~\ref{7.3} is rather difficult and requires some preliminary work.

A base $\mathcal B$ of the topology of a space $X$ will be called {\em regular} if for each open set $U\subset X$ and a point $x\in U$ there are two sets $V,W\in\mathcal B$ with $x\in V\subset X\setminus W\subset U$. Such a regular base $\mathcal N$ will be called {\em countably additive} if the union $\cup\C$ of any countable subfamily $\C\subset\mathcal B$ belongs to $\mathcal B$. 

Observe that a topological space $X$ admits a regular base if and only if $X$ is regular.
Note also that the family of all functionally open subsets of a Tychonov space $X$ forms a regular countably additive base of the topology of $X$. We recall that a subset $U\subset X$ is called {\em functionally open} if $U=f^{-1}(V)$ for some continuous map $f:X\to\IR$ and some open set $V\subset\IR$. 

We define a map $f:X\to Y$ to be {\em weakly $G_\delta$-measurable} if there is a regular countably additive base $\mathcal B$ of the topology of $Y$ such that for every $U\in\mathcal B$ the preimage $f^{-1}(U)$ is of type $G_\delta$ in $X$. 

Since the topology of each regular space forms a regular countably additive base, we conclude that each $G_\delta$-measurable map $f:X\to Y$ into a regular topological space $Y$ is weakly $G_\delta$-measurable.

Let us say that a map $f:X\to Y$ between topological spaces is
{\em almost continuous} ({\em quasi-continuous}) at a point $x\in X$ if for any neighborhood $Oy\subset Y$ of the point $y=f(x)$ the (interior of the) set $f^{-1}(Oy)$ is dense in some neighborhood of the point $x$ in $X$. By $AC(f)$ (resp. $QC(f)$~) we shall denote the set of point of almost (resp.  quasi-) continuity of $f$.

\begin{lem}\label{7.2} Let $f:X\to Y$ be a weakly $G_\delta$-measurable map from a Baire space $X$ to a regular  space $Y$. Then
\begin{enumerate}
\item $AC(f)=QC(f)$.
\item If $D$ is dense in $X$ and $f|D$ has no continuity point, then $D\setminus AC(f)$ also is dense in $X$.
\item For any countable dense set $D\subset X$ there is a point $y\in f(D)$ such that for every neighborhood $Oy$ of $y$ the preimage $f^{-1}(Oy)$ has non-empty interior in $X$.
\item The family $\{\Int \overline{f^{-1}(y)}:y\in Y\}$ is disjoint.
\end{enumerate}
\end{lem}

\begin{proof} Let $\mathcal B$ be a regular countably additive base of the topology of $Y$ such that for every $U\in\mathcal B$ the preimage $f^{-1}(U)$ is of type $G_\delta$ in $X$.

1. The inclusion $QC(f)\subset AC(f)$ is trivial. To prove that $AC(f)\subset QC(f)$, take any point $x\in AC(f)$. To show that $x\in QC(f)$, take any neighborhood $Oy\subset Y$ of the point $y=f(x)$. 
We should check that the interior of $f^{-1}(Oy)$ is dense in some neighborhood of $x$.
By the regularity of the base $\mathcal B$, there are sets $Vy,Wy\in\mathcal B$ such that $y\in Vy\subset Y\setminus Wy\subset Oy$. Then the preimages $f^{-1}(Vy)$ and $f^{-1}(W_y)$ are of type $G_\delta$ and $F_\sigma$ in $X$, respectively. Since $x\in AC(f)$, the closure of $f^{-1}(Vy)$ contains an open neighborhood $Ox$ of $x$. We claim that the closure of the interior of $f^{-1}(Oy)$ contains $Ox$. This will follow as soon as, given a non-empty open subset  $U\subset Ox$ we find a non-empty open subset $U'\subset U\cap f^{-1}(Oy)$. Observe that the set $U\cap f^{-1}(Vy)$, being a dense $G_\delta$-set in the Baire space $U$, is not meager in $U$. Consquently, the $F_\sigma$-subset $U\cap f^{-1}(Wy)\supset O_x\cap f^{-1}(Vy)$ is non-meager in $U$ too. Now the Baire Theorem yields a non-empty open subset $$U'\subset U\cap f^{-1}(Wy)\subset U\cap f^{-1}(Oy).$$
\smallskip

2. Assume that $D\subset X$ is dense and $f|D$ has no continuity point. Given a point $x\in D$, and a neighborhood $O_x$ of $x$ we should find a point $x'\in O_x\cap D\setminus AC(f)$. If $x\notin AC(f)$, then we can take $x'=x$. So we assume that $x\in AC(f)$ and hence $x\in QC(f)$ by the preceding item. Since $x$ is a discontinuity point of $f|D$, there is a neighborhood $O_{f(x)}$ of $f(x)$ such that $f(D\cap U_x)\not\subset O_{f(x)}$ for every neighborhood $U_x$ of $x$. Using the regularity of $Y$ choose a neighborhood $U_{f(x)}\subset Y$ of $f(x)$ with $\overline{U_{f(x)}}\subset O_{f(x)}$. Since $x$ is a quasi-continuity point of $f$, the closure  of the interior  of the preimage $f^{-1}(U_{f(x)})$ contains some neighborhood $W_x$ of $x$. By the choice of $O_{f(x)}$, we can find a point $x'\in D\cap O_x\cap W_x$ with $f(x')\notin O_{f(x)}$. Consider the neighborhood $O_{f(x')}=Y\setminus \overline{U_{f(x)}}$ of $f(x')$ and observe that $W_x\cap f^{-1}(O_{f(x')})$ is a nowhere dense subset of $Ox$ (because it misses the interior of $f^{-1}(U_{f(x)})$ which is dense in $W_x$). This witnesses that $x'\notin AC(f)$.
\smallskip

3. Given a countable dense subset $D\subset X$, we should find a point $y\in f(D)$ such that for every neighborhood $Oy$ the preimage $f^{-1}(Oy)$ has non-empty interior in $X$.

Assume conversely that each point $y\in f(D)$ has a neighborhood $Oy$ such that the preimage $f^{-1}(Oy)$ has empty interior in $X$. Using the regularity of $\mathcal B$, for every $y\in f(D)$ find two sets $Vy,Wy\in\mathcal B$ with $y\in Vy\subset Y\setminus Wy\subset Oy$. It follows that $f^{-1}(Vy)$ is of type $G_\delta$ while $f^{-1}(Y\setminus Wy)$ of type $F_\sigma$ in $X$.
Being an $F_\sigma$-set with empty interior in the Baire space $X$, the set $f^{-1}(Y\setminus Wy)$ is meager in $X$. Then the union $$U=\bigcup_{y\in f(D)}f^{-1}(Vy)\supset D$$ is a dense meager set in $X$ too. 

On the other hand, the set $U$ is a $G_\delta$-subset of $X$, being the preimage $U=f^{-1}(V)$ of the subset $V=\bigcup_{y\in f(D)}Vy$ that belongs to $\mathcal B$ because of the countable additivity of $\mathcal B$. By the Baire Theorem, the dense $G_\delta$-subset $U$ of the Baire space $X$ cannot be meager, which is a contradiction.
\smallskip

4. Assuming that the family $\{\Int \overline{f^{-1}(y)}:y\in Y\}$ is not disjoint, find two distinct points $y,z\in Y$ such that the intersection  
$$W=\Int \overline{f^{-1}(y)}\cap \Int \overline{f^{-1}(z)}$$is not empty.
Observe that the sets $W\cap f^{-1}(y)$ and $W\cap f^{-1}(z)$ both are dense in $W$.

By the Hausdorff property of $Y$ the points $y,z$ have disjoint open neighborhoods $Oy,Oz\in\mathcal B$.  The choice of $\mathcal B$ guarantees that the sets $W\cap f^{-1}(Oy)$ and $W\cap f^{-1}(Oz)$ are dense disjoint $G_\delta$-sets in the Baire space $W$, which is forbidden by the Baire Theorem. 
\end{proof}

\begin{lem}\label{7.2a} Let $f:X\to Y$ be a weakly $G_\delta$-measurable map from a Baire space $X$ to a regular  space $Y$ and $D$ be a countable dense subset of $X$ such that $f|D$ has no continuity point.
\begin{enumerate}
\item For any finite subset $F\subset Y$ there is a dense subset $Q\subset D\setminus f^{-1}(F)$ in $X$ such that $f|Q$ has no continuity point.
\item For any sequence $(U_n)_{n=1}^\infty$ of non-empty open subsets of $X$ there are an infinite subset $I\subset\IN$ and sequences $(V_n)_{n\in I}$ and $(W_n)_{n\in I}$ of pairwise disjoint non-empty open sets in $X$ and $Y$, respectively, such that  $V_n\subset U_n\cap f^{-1}(W_n)$ for all $n\in I$.
\item If $D\subset PS(X)$, then there is a countable first countable subspace $Q\subset D$ having no isolated points and such that the restriction $f|Q$ is a bijective map whose image $f(Q)$ is a discrete subspace of $Y$.
\end{enumerate}
\end{lem}

\begin{proof} Fix a regular countably additive base $\mathcal B$ of the topology of $Y$ such that for every $U\in\mathcal B$ the preimage $f^{-1}(U)$ is of type $G_\delta$ in $X$.
\smallskip

1. Take any finite subset $F=\{y_1,\dots,y_n\}$ in $Y$ and for every $i\le n$ let $D_i=D\cap f^{-1}(y_i)$ and $W_i=\Int \overline{D_i}$. 

First we show that
the complement $D\setminus f^{-1}(F)$ is dense in $X$. Assuming the converse, we can find a non-empty open set $U\subset X$ such that $U\cap D\subset f^{-1}(F)$. 
We claim that $U\cap W_i\cap W_j\ne\emptyset$ for some distinct numbers $i,j\le n$. Assuming conversely that the family $\{U\cap W_i\}_{i=1}^n$ is disjoint, we would conclude that $W_1\cap U\cap D\subset D_1$ and hence the restriction $f|D\cap W_1\cap U$, being a constant map to $y_1$, is continuous. But this contradicts the fact that $f|D$ has no continuity point. So, $W_0=U\cap W_i\cap W_j$ is not empty for some distinct $i,j\le n$. Then $D_i\cap W_0$ and $D_j\cap W_0$ are disjoint dense subsets in $W_0$. Using the Hausdorff property fo $Y$, take two disjoint neighborhoods $Oy_i,Oy_j\in\mathcal B$. Then the dense disjoint sets $G_i=W_0\cap f^{-1}(Oy_i)\supset W_0\cap D_i$ and  
$G_j=W_0\cap f^{-1}(Oy_j)\supset W_0\cap D_j$ are of type $G_\delta$ in the Baire space $W_0$, which contradicts the Baire Theorem. This completes the proof of the density of $D\setminus f^{-1}(F)$.

It follows that the set
$$Q=\big(D\setminus f^{-1}(F)\big)\setminus\bigcup_{i=1}^n \overline{W_i}\setminus W_i$$
is  dense in $X$. We claim that the restriction $f|Q$ has no continuity point.
Assume conversely that some point $x_0\in Q$ is a continuity point of the restriction $f|Q$. If $x_0\notin\bigcup_{i=1}^n\overline{W_i}$, then the discontinuity of the map $f|D$ at $x_0$ implies the discontinuity of $f|Q$ at $x_0$. So, $x_0\in W_i$ for some $i\le n$. Let $y_0=f(x_0)$ and observe that $y_0\ne y_i$ (because $x_0\notin f^{-1}(F)$). By the Hausdorff property of $Y$ the points $y_0$ and $y_i$ have disjoint open neighborhoods $Oy_0,Oy_i\in\mathcal B$.
By the continuity of $f|Q$ at $x_0$, there is an open neighborhood $Ox_0\subset W_i$ of $x_0$ such that $f(Ox_0\cap Q)\subset Oy_0$. It follows that 
 the sets $G_0=Ox_0\cap f^{-1}(Oy_0)$ and $G_i=Ox_0\cap f^{-1}(Oy_i)$ are of type $G_\delta$ in $Ox_0$. The density of the set $D\setminus f^{-1}(F)$ implies the density of the set
$$G_0\supset Ox_0\cap Q=(D\setminus f^{-1}(F))\cap W_i$$ in $Ox_0$.
On the other hand, the definition of the set $W_i$ implies that the intersection $f^{-1}(y_i)\cap W_i$ is dense in $W_i$ and hence $$G_i\supset f^{-1}(y_i)\cap Ox_0$$ is dense in $Ox_0$. Thus we have found two disjoint dense $G_\delta$-sets in the Baire space $Ox_0$ which is forbidden by the Baire Theorem.
\smallskip

2. Let $(U_n)_{n=1}^\infty$ be a sequence of non-empty open subsets of $X$.
Applying Lemma~\ref{7.2}(3) to the map $f|U_1$ and the dense subset $D\cap U_1$, find a point $y_0\in f(D\cap U_1)$ such that for each neighborhood $Oy_0$ the preimage $U_1\cap f^{-1}(Oy_0)$ has non-empty interior. 
By induction, for every $n\in\IN$ we shall find a point $y_n\in f(D)\setminus \{y_i:i<n\}$ such that for every neighborhood $Oy_n$ the set $U_n\cap f^{-1}(Oy_n)$ has non-empty interior. 

Assuming that for some $n$ the points $y_0,\dots,y_{n-1}$ have beeing chosen, we shall find a point $y_n$. It follows that the intersection $D\cap U_n$ is a countable dense subset of $U_n$ such that $f|D\cap U_n$ has no continuity point. Applying the preceding item, we can find a dense subset $Q\subset D\cap U_n\setminus f^{-1}(\{y_0,\dots,y_{n-1}\})$ in $U_n$ such that the restriction $f|Q$ has no continuity point. Applying Lemma~\ref{7.2}(3) to the map $f|U_n$ and the dense subset $Q$, find a point $y_n\in f(Q)\subset f(D)\setminus\{y_i:i<n\}$ such that for each neighborhood $Oy_n$ the preimage $U_n\cap f^{-1}(Oy_n)$ has non-empty interior. 
This completes the inductive construction.
\smallskip

The space $\{y_n:n\in\IN\}$, being infinite and regular, contains an infinite discrete subspace $\{y_n:n\in I\}$. By induction, we can select pairwise disjoint open neighborhoods $W_n\subset Y$, $n\in I$, of the points $y_n$. For every $n\in I$ the set $U_n\cap f^{-1}(W_n)$ contains a non-empty open set $V_n$ by the choice of the point $y_n$. Then the set $I\subset\IN$ and sequences $(V_n)_{n\in I}$, $(W_n)_{n\in I}$ satisfy our requirements.
\smallskip

3. Assume that $D\subset PS(X)$. Applying Lemma~\ref{7.2}(2), we get that $D\setminus AC(f)$ is dense in $X$.

By induction on the tree $\w^{<\w}$ we shall construct sequences $(x_s)_{s\in\w^{<\w}}$ of points of the set $D\setminus AC(f)$, and sequences $(V_s)_{s\in \w^{<\w}}$ and $(U_s)_{s\in\w^{<\w}}$, $(W_s)_{s\in\w^{<\w}}$ of sets so that the following conditions hold for every finite number sequence $s\in\w^{<\w}$:
\begin{itemize}
\item[(a)] $V_s$ is an open neighborhood of the point $x_s$ in $X$;
\item[(b)] $W_s\subset U_s$ are open neighborhoods of $f(x_s)$ in $Y$;
\item[(c)] $f(V_s)\subset U_s$;
\item[(d)] $V_{s\concat n}\subset V_s$ and $U_{s\concat n}\subset U_s$ for all $n\in\w$;
\item[(e)] the sequence $(V_{s\concat n})_{n\in\IN}$ converges to $x_s$;
\item[(f)] $W_{s}\cap U_{s\concat n}=\emptyset=U_{s\concat n}\cap U_{s\concat m}$ for all $n\ne m$ in $\w$.
\end{itemize}

We start the induction letting $V_\emptyset=X$, $U_\emptyset=Y$ and $x_\emptyset$ be any point of $D\setminus AC(f)$.

Assume that for a finite sequence $s\in\w^{<\w}$ the point $x_s\in D\setminus AC(f)$
and open sets $V_s\subset X$ and $U_s\subset Y$ with $x_s\in V_s$ and $f(V_s)\subset U_s$ have been constructed.  Since $f|V_s$ fails to be almost continuous at $x_s$, there is a neighborhood $W_s\subset U_s$ of $f(x_s)$ such that the closure of the preimage $f^{-1}(\overline{W_s})$ is not a neighborhood of $x_s$ in $X$. This fact and the Preiss-Simon property of $X$ at $x_s$ allows us to construct a sequence $(V'_{k})_{k\in\w}$ of open subsets of $V_s\setminus \cl_X\big(f^{-1}(\overline{W_s})\big)$ that converges to $x_s$ in the sense that each neighborhood of $x$ contains all but finitely many sets $V'_{k}$. Applying item 2 to the map $f|V_s:V_s\to U_s$,
we can find an infinite subset $N\subset \w$ and a sequence $(U'_{k})_{k\in N}$ of pairwise disjoint open sets of $U_s$ such that each set $f^{-1}(U'_k)\cap V'_k$, $k\in N$, has non-empty interior in $X$. Let $N=\{k_n:n\in\w\}$ be the increasing enumeration of the set $N$.

For every $n\in\w$ let $U_{s \concat n}=U'_{k_n}\setminus\overline{W_s}$, $V_{s\concat n}$ be a non-empty open subset in $f^{-1}(U_{k_n})\cap V_{k_n}'$ and $x_{s\concat n}\in V_{s\concat n}\cap D\setminus AC(f)$ be any point (such a point exists because of the density of $D\setminus AC(f)$ in $X$). One can check that the points $x_{s\concat n}$, ${n\in\w}$ and sets $W_s$, $V_{s\concat n}$, $U_{s\concat n}$, $n\in\w$ satisfy the requirements of the inductive construction.

After completing the inductive construction, consider the set $Q=\{x_s:s\in\w^{<\w}\}$ and note that it is first countable and has no isolated point, $f|Q$ is bijective and $f(Q)$ is a discrete subspace of $Y$.
\end{proof}

Now we are able to prove the promised

\begin{thm}\label{7.3} Each weakly $G_\delta$-measurable map $f:X\to Y$ from a hereditarily Baire Preiss-Simon space $X$ to a regular space $Y$ is scatteredly continuous.
\end{thm}

\begin{proof} Fix a countably additive regular base $\mathcal B$ of the topology of $Y$ such that $f^{-1}(U)$ is of type $G_\delta$ in $X$ for every $U\in\mathcal B$. Assuming that $f$ is not scatteredly continuous we could find a subspace $D\subset X$ such that the restriction $f|D$ has no continuity point. Being Preiss-Simon, the space $X$ has countable tightness. In this case we can additionally assume that the set $D$ is countable, see Proposition~\ref{2.4}.  Applying
Lemma~\ref{7.2a}(3) to the restriction $f|\overline{D}$ of $f$ onto the Baire space $\overline{D}$, we can find a countable subset $Q\subset D$ without isolated point such that $f|Q$ is bijective and $f(Q)$ is a discrete subspace of $Y$. It is clear that $f|Q$ has no continuity point.
Observe that $Z=\overline{Q}$ is a Baire space, being a closed subspace of the hereditarily Baire space $X$. Let $g=f|Z$ and observe that for every $U\in\mathcal B$ the preimage $g^{-1}(U)$ is of type $G_\delta$ in $Z$.

Since the set $g(Q)=f(Q)$ is countable and discrete, for every $y\in g(Q)$ we can select a neighborhood $Oy\in\mathcal B$ so small that the family $\{Oy:y\in g(Q)\}$ is disjoint. The countable additivity of the base $\mathcal B$ implies that $\bigcup_{y\in g(Q)}Oy\in\mathcal B$ and hence the
preimage $G=g^{-1}(\bigcup_{y\in g(Q)}Oy)$ is a Baire space, being a dense $G_\delta$-subset of the Baire space $Z$. Since $G=\bigcup_{y\in g(Q)}g^{-1}(Oy)$, the Baire Theorem guarantees that for some $z\in g(Q)$ the preimage $g^{-1}(Oz)$ is not meager in $G$. On the other hand, $\bigcup_{y\in g(Q)\setminus\{z\}}Oy\in\mathcal B$ and hence the complement $G\setminus g^{-1}(Oz)=g^{-1}(\bigcup_{y\in g(Q)\setminus\{z\}}Oy)$ is a $G_\delta$-set in $G$.
The set $g^{-1}(Oz)$, being a non-meager $F_\sigma$-subset of the Baire space $G$, contains a non-empty open subset $W$ of $G$. For this set we get $$g(W\cap Q)\subset Oz\cap g(Q)=\{z\},$$ which means that $g|Q$ is continuous at points of $W\cap Q$. But this contradicts the fact that $g|Q\equiv f|Q$ has no continuity point.
\end{proof}

\section{Characterizing the scattered continuity}

Unifying Theorems~\ref{4.4}, \ref{5.1}, \ref{6.1}, \ref{7.3} and Propositions~\ref{6.0}, \ref{7.1}, we obtain the following characterization of scatteredly continuous maps:

\begin{thm}\label{8.1} For a map $f:X\to Y$ from a perfectly paracompact hereditarily Baire  space $X$ to a regular space $Y$ the following conditions are equivalent:
\begin{enumerate}
\item $f$ is scatteredly continuous;
\item $f$ is weakly discontinuous;
\item $f$ is piecewise continuous.
\end{enumerate}
If $X$ is a Preiss-Simon space, then (1)--(3) are equivalent to
\begin{enumerate}
\item[(4)] $f$ is $G_\delta$-measurable;
\item[(5)] $f$ is weakly $G_\delta$-measurable;
\end{enumerate}
If $Y$ is path-connected and belongs to $\sAE(X)$ , then (1)--(3) are equivalent to
\begin{enumerate}
\item[(6)] $f$ is of stable first Baire class.
\end{enumerate}
\end{thm}

Theorem~\ref{8.1} has an interesting corollary describing the interplay between $F_\sigma$- and $G_\delta$-measurable maps. First observe that, for a space $Y$ whose each open subset is of type $F_\sigma$, each $G_\delta$-measurable map $f:X\to Y$ is $F_\sigma$-measurable. The same effect can be achieved imposing some restrictions on the domain $X$ of $f$. 

\begin{cor} Each (weakly) $G_\delta$-measurable map $f:X\to Y$ from a perfectly paracompact hereditarily Baire Preiss-Simon space $X$ to a regular space $Y$ is  $F_\sigma$-measurable.
\end{cor}

\begin{question} Let $f:X\to Y$ be a $G_\delta$-measurable function from a perfectly paracompact hereditarily Baire space $X$ to a regular space $Y$. Is $f$ $F_\sigma$-measurable?
\end{question}  

\section{Some Examples}\label{s9}

The requirement that $X$ is hereditarily Baire is necessary in Theorem~\ref{8.1}.

\begin{prop} If a perfectly paracompact space $X$ is not hereditarily Baire, then there is a bijective piecewise continuous map $f:X\to Y$ that is not scatteredly continuous.
\end{prop}

\begin{proof} Since the space $X$ is not hereditarily Baire, there is a closed subspace $F\subset X$ of the first Baire category in itself. Write $F=\bigcup_{n\in\w}F_n$ as the union of an increasing sequence $(F_n)_{n\in\w}$ of closed nowhere dense subsets of $F$ with $F_0=\emptyset$. Let $Y$ be the topological sum of the family $\F=\{X\setminus F,F_{n+1}\setminus F_n:n\in\w\}$ and $f:X\to Y$ be the identity map. The perfect paracompactness of $X$ yields that each subset $A\in \F$ is of type $F_\sigma$ in $X$, which implies  that $f$ is piecewise continuous. On the other hand, the restriction $f|F$ has no continuity point, which means that $f$ fails to be scatteredly continuous.
\end{proof}

The above proposition shows that the equivalences $(1)\Leftrightarrow(2)\Leftrightarrow(3)$ of Theorem~\ref{8.1} do not hold without the assumption that $X$ is hereditarily Baire. On the other hand, the equivalence $(3)\Leftrightarrow(4)$ does not require that assumption and hold, for example, for any map $f:X\to Y$ from an absolute $\F$-Souslin space $X$ to a metrizable space $Y$, see \cite{JR82}. 

Now we see that a $G_\delta$-measurable function $f:X\to Y$ between metrizable spaces is piecewise continuous provided $X$ is either hereditarily Baire (by Theorem~\ref{7.3}) or absolute $\F$-Souslin (by \cite{JR82}). In this situation it is natural to ask if any $G_\delta$-measurable function $f:X\to Y$ between metrizable spaces is piecewise continuous, see V.Vinokurov \cite{Vino}. Below we give two (consistent) counterexamples to this question.

The first one uses $Q$-spaces. Those are topological spaces $X$ such that each subset $A\subset X$ is of type $F_\sigma$ in $X$. The Martin Axiom implies that each second countable space $X$ of cardinality 
$|X|<\mathfrak c$ is a $Q$-space, see \cite[4.2]{Miller}. On the other hand, in some models of ZFC, each second countable $Q$-space is at most countable, see \cite[4.3]{Miller}.

\begin{ex}\label{9.2} If $X$ is an uncountable second countable $Q$-space, then any bijective map $f:X\to D$ to a discrete space $D$ is $G_\delta$-measurable but not piecewise continuous. Moreover, $\dec(f)=|X|$.
\end{ex}

Here for a function $f:X\to Y$ by $\dec(f)$ we denote the smallest cardinality $|\C|$ of a cover $\C$ of $X$ such that $f|C$ is continuous for every $C\in\C$, see \cite{Sol}. It is clear that $\dec(f)\le\dec_c(f)$. The latter inequality can be strict: for the characteristic function $\chi_\IQ:\IR\to\{0,1\}$ of rationals we get $\dec(f)=\aleph_0$ while $\dec_c(f)=\mathfrak d$, see \cite{Sol}.

A typical example of an $F_\sigma$-measurable function with $\dec(f)>\aleph_0$ is the Pawlikowski function $P:(\w+1)^\w\to \w^\w$, which is the countable power of the bijection $\w+1\to\w$ sending every $n\in\w$ to $n+1$ and $\w$ to $0$, see \cite{CMPS}, \cite{Sol}. Here the ordinal $\w+1=\w\cup\{\w\}$, endowed with the order topology, is homeomorphic to a convergent sequence.

Note that the range space $D$ in Example~\ref{9.2} is not second countable. Our second counterexample to the Vinokurov question has second countable range. It relies on the notion of a Lusin space.

Given a cardinal $\kappa$ we define a topological space $X$ to be {\em $\kappa$-Lusin} if each nowhere dense subset $A$ of $X$ has size $|A|<\kappa$. Uncountable $\aleph_1$-Lusin spaces are referred to as Lusin spaces, see \cite[\S2]{Miller}. We shall be interested in $\mathfrak b$-Lusin spaces, where  
 $\mathfrak b$ is the smallest cardinality $|A|$ of a subset $A\subset\IN^\w$ that does not lie in a $\sigma$-compact subset of $\IN^\w$, see \cite[\S2]{Blass}.

\begin{ex}\label{9.3} Let $P:(\w+1)^\w\to\w^\w$ be the Pawlikowski function. For every $\bb$-Lusin subset $L\subset\w^\w$ of size $|L|\ge\bb$ and its preimage $X=P^{-1}(L)$ the function $f=P|X:X\to \w^\w$ is $G_\delta$-measurable but has $\dec(f)=|L|\ge\bb$.
\end{ex}

\begin{proof} The $G_\delta$-measurability of $f=P|X$ will follow as soon as we check that for every closed subset $F\subset\w^\w$ the set $f^{-1}(F)=P^{-1}(F)\cap X$ is of type $F_\sigma$ in $X$. First we consider the case of a nowhere dense $F$. In this case $|F\cap L|<\bb$ and hence the set $A=f^{-1}(F)=P^{-1}(F)\cap X$ has size $|A|=|L\cap F|<\bb$. The space $P^{-1}(F)$, being a Borel subset of $(\w+1)^\w$ is analytic and thus is the image of $\w^\w$ under a suitable continuous map $g:\w^\w\to P^{-1}(F)$. Find a subset $A'\subset\w^\w$ of size $|A'|=|A|$ with $g(A')=A$. The definition of the cardinal $\bb$ implies the existence of a $\sigma$-compact subset $K\subset\w^\w$ with $K\supset A'$. Then $g(K)$ is a $\sigma$-compact subset of $P^{-1}(F)$ containing the set $A$. Consequently, $A=P^{-1}(F)\cap X=g(K)\cap X$ is of type $F_\sigma$ in $X$. 

If $F$ is an arbitrary closed subset of $\w^\w$, then take the interior $U$ of $F$ in $\w^\w$ and observe that $P^{-1}(U)$ is of type $F_\sigma$ in $(\w+1)^\w$ by the $F_\sigma$-measurability of $P$. Then $f^{-1}(U)=P^{-1}(U)\cap X$ is of type $F_\sigma$ in $X$. Since $F\setminus U$ is closed and  nowhere dense in $\w^\w$, the preimage $f^{-1}(F\setminus U)$ is of type $F_\sigma$ in $X$ too. Combining these two facts, we conclude that the preimage $f^{-1}(F)=f^{-1}(U)\cup f^{-1}(F\setminus U)$ is of type $F_\sigma$ in $X$, witnessing that $f:X\to \w^\w$ is $G_\delta$-measurable.

Next, we check that $\dec(f)\ge|L|$. Let $\C$ be a cover of $X$ such that $|\C|=\dec(f)$ and $f|C$ is continuous for every $C\in\C$. Using the definition of the Pawlikowski function, it is easy to check that the image $f(C)$ of each $C\in\C$ is nowhere dense in $\w^\w$ and hence $|f(C)\cap L|<\bb$. It follows from $|L|\ge\bb$ and the regularity of the cardinal $\bb$ \cite[2.4]{Blass} that $|\C|\ge|L|$.
\end{proof}

In light of Example~\ref{9.3} it is important to have some conditions guaranteeing the existence of $\bb$-Lusin sets of cardinality $\ge\bb$.
Such conditions can be written in terms of the following three cardinal characteristics of the  $\sigma$-ideal $\M$ of meager subsets of the Baire space $\w^\w$:
$$
\begin{aligned}
\add(\M)=&\min\{|\A|:\A\subset\M,\;\cup\A\notin\M\},\\
\cov(\M)=&\min\{|\A|:\A\subset\M,\; \cup\A=\w^\w\},\\
\cof(\M)=&\min\{|\A|:\A\subset\M\;\;\forall M\in\M\;\exists A\in\A\mbox{ with }M\subset A\}
\end{aligned}
$$

The following proposition is a modification of Corollary~8.2.5 and Lemma~8.2.4 from \cite{BJ}.

\begin{prop} \begin{enumerate}
\item If $\add(\M)=\cof(\M)$, then there is a $\bb$-Lusin subset $L\subset\w^\w$ of size $|L|=\bb$;
\item If there is a $\bb$-Lusin subset $L\subset\w^\w$ of size $|L|\ge\bb$, then $\cov(\M)\ge\bb$.
\end{enumerate}
\end{prop}

\begin{proof} 1. The Cicho\'n diagram \cite[\S5]{Blass} guarantees $\add(\M)\le\bb\le\cof(\M)$.   Consequently, the equality $\add(\M)=\cof(\M)$ implies $\add(\M)=\bb=\cof(\M)$. Let $\{M_\alpha\}_{\alpha<\bb}\subset\M$ be a cofinal family of meager subsets of $\w^\w$. By the transfinite induction, for every $\alpha<\bb$ choose a point $x_\alpha\notin \bigcup_{\beta\le \alpha
}M_\beta$ in $\w^\w$. Such a choice is always possible, since the latter union is meager (by $\add(\M)=\bb$). The set $L=\{x_\alpha:\alpha<\bb\}$ is the required $\bb$-Lusin subset of $\w^\w$ with $|L|=\bb$.
\smallskip

2. Now assume that $L\subset\w^\w$ is a $\bb$-Lusin subset with $|L|\ge \bb$. Fix a cover $\C\subset\M$ of $\w^\w$ of size $|C|=\cov(\M)$. It follows that $|L\cap C|<\bb$ for all $C\in\C$. Now the regularity of the cardinal $\bb$ \cite[2.4]{Blass} implies that $\cov(\M)=|\C|\ge\bb$.
\end{proof}

By \cite[7.4]{Blass}, the Martin Axiom implies $\add(\M)=\cof(\M)=\cc$. On the other hand, the strict inequality $\cov(\M)<\bb$ holds is the Laver's model of ZFC, see \cite[11.7]{Blass}. In this model no $\bb$-Lusin subset $L\subset\w^\w$ with $|L|\ge\bb$ exists. This shows that Example~\ref{9.3} has consistent nature.

\begin{question} Is it consistent with ZFC that each $G_\delta$-measurable function $f:X\to\IR$ defined on a subset $X\subset\IR$ has $\dec(f)\le\aleph_0$? $\dec_c(f)\le\aleph_0$?
\end{question} 

\section{Acknowledgement}

The authors express their sincere thanks to Lubomyr Zdomskyy for fruitful discussions resulting in finding Example~\ref{9.3}.


\end{document}